\documentclass[12pt]{amsart}
\usepackage{amssymb}

\newcommand{\Bgp}{{\Z^\N}}

\long\def\forget#1\forgotten{}
\newcommand{\issuenumber}{28}
\newcommand{\issuemonth}{August}
\newcommand{\issueyear}{2009}

\newcommand{\Fin}{{[\N]^{<\aleph_0}}}
\newcommand{\alephes}{{\aleph_0}}
\newcommand{\ed}{
\newpage

\section{Unsolved problems from earlier issues}

\begin{issue}
Is $\binom{\Omega}{\Gamma}=\binom{\Omega}{\Tau}$?
\end{issue}

\begin{issue}
Is $\ufin(\cO,\Omega)=\sfin(\Gamma,\Omega)$?
And if not, does $\ufin(\cO,\Gamma)$ imply
$\sfin(\Gamma,\Omega)$?
\end{issue}

\stepcounter{issue}

\begin{issue}
Does $\sone(\Omega,\Tau)$ imply $\ufin(\Gamma,\Gamma)$?
\end{issue}

\begin{issue}
Is $\fp=\fp^*$? (See the definition of $\fp^*$ in that issue.)
\end{issue}

\begin{issue}
Does there exist (in ZFC) an uncountable set satisfying $\sfin(\B,\B)$?
\end{issue}

\stepcounter{issue}

\begin{issue}
Does $X \nin \NON(\cM)$ and $Y\nin\mathsf{D}$ imply that
$X\cup Y\nin \COF(\cM)$?
\end{issue}

\begin{issue}[CH]
Is $\split(\Lambda,\Lambda)$ preserved under finite unions?
\end{issue}

\begin{issue}
Is $\cov(\cM)=\fo$? (See the definition of $\fo$ in that issue.)
\end{issue}

\begin{issue}
Does $\sone(\Gamma,\Gamma)$ always contain an element of cardinality $\fb$?
\end{issue}

\begin{issue}
Could there be a Baire metric space $M$ of weight $\aleph_1$ and a partition
$\mathcal{U}$ of $M$ into $\aleph_1$ meager sets where for each ${\mathcal U}'\subset\mathcal U$,
$\bigcup {\mathcal U}'$ has the Baire property in $M$?
\end{issue}

\stepcounter{issue} 

\begin{issue}
Does there exist (in ZFC) a set of reals $X$ of cardinality $\fd$ such that all
finite powers of $X$ have Menger's property $\sfin(\cO,\cO)$?
\end{issue}

\begin{issue}
Can a Borel non-$\sigma$-compact group be generated by a Hurewicz subspace?
\end{issue}

\begin{issue}[MA]
Is there $X\sbst\R$ of cardinality continuum, satisfying $\sone(\BO,\BG)$?
\end{issue}

\begin{issue}[CH]
Is there a totally imperfect $X$ satisfying $\ufin(\cO,\Gamma)$
that can be mapped continuously onto $\Cantor$?
\end{issue}

\begin{issue}[CH]
Is there a Hurewicz $X$ such that $X^2$ is Menger but not Hurewicz?
\end{issue}

\begin{issue}
Does the Pytkeev property of $C_p(X)$ imply that $X$ has Menger's property?
\end{issue}

\begin{issue}
Does every hereditarily Hurewicz space satisfy $\sone(\BG,\BG)$?
\end{issue}

\begin{issue}[CH]
Is there a Rothberger-bounded $G\le\Bgp$ such that $G^2$ is not Menger-bounded?
\end{issue}

\begin{issue}
Let $\cW$ be the van der Waerden ideal.
Are $\cW$-ultrafilters closed under products?
\end{issue}

\begin{issue}
Is the $\delta$-property equivalent to the $\gamma$-property $\binom{\Omega}{\Gamma}$?
\end{issue}

\stepcounter{issue}

\stepcounter{issue}

\general\end{document}}

\newcommand{\Cantor}{{\{0,1\}^\N}}

\newcommand{\roth}{{[\N]^{\alephes}}}

\newcommand{\fb}{\mathfrak{b}}

\newcommand{\fc}{\mathfrak{c}}
\newcommand{\fd}{\mathfrak{d}}
\newcommand{\fp}{\mathfrak{p}}

\newcommand{\NON}{{\mathsf   {NON}}}
\newcommand{\COF}{{\mathsf   {COF}}}

\newcommand{\cM}{\mathcal{M}}

\newcommand{\cov}{\mathsf{cov}}

\newcommand{\non}{\mathsf{non}}

\newcommand{\CH}{the Continuum Hypothesis}
\newcommand{\R}{\mathbb{R}}

\newcommand{\fo}{\mathfrak{od}}

\newcommand{\ft}{\mathfrak{t}}

\renewcommand{\split}{\mathsf{Split}}
\newcommand{\bq}{\begin{quote}}
\newcommand{\eq}{\end{quote}}
\newcommand{\cO}{\mathcal{O}}
\newcommand{\B}{\mathcal{B}}
\newcommand{\BG}{\B_\Gamma}

\newcommand{\BO}{\B_\Omega}

\newcommand{\sone}{\mathsf{S}_1}    \newcommand{\sfin}{\mathsf{S}_\mathrm{fin}}

\newcommand{\ufin}{\mathsf{U}_\mathrm{fin}}

\newcommand{\nin}{\not\in}
\newcommand{\cF}{\mathcal{F}}

\newcommand{\cU}{\mathcal{U}}
\newcommand{\cV}{\mathcal{V}}
\newcommand{\cW}{\mathcal{W}}

\newcommand{\NN}{{\N^\N}}
\newcommand{\N}{\mathbb{N}}
\newcommand{\Z}{\mathbb{Z}}
\newcommand{\as}{\subseteq^*}
\newcommand{\sm}{\setminus}
\newcommand{\sbst}{\subseteq}
\newcommand{\by}[2]{\par\hfill\emph{#1}, #2}
\newcommand{\nby}[1]{\par\hfill\emph{#1}}
\newcommand{\Tau}{\mathrm{T}}
\newcommand{\CE}{\textsc{CE}}

\newtheorem{thm}{Theorem}[section]
\newcommand{\bthm}{\begin{thm}} \newcommand{\ethm}{\end{thm}}
\newtheorem{prop}[thm]{Proposition}
\newcommand{\bprp}{\begin{prop}} \newcommand{\eprp}{\end{prop}}
\newtheorem{fact}[thm]{Fact}
\newcommand{\bfct}{\begin{fact}} \newcommand{\efct}{\end{fact}}
\newtheorem{prob}[thm]{Problem}
\newcommand{\bprb}{\begin{prob}} \newcommand{\eprb}{\end{prob}}
\newtheorem{lem}[thm]{Lemma}
\newcommand{\blem}{\begin{lem}} \newcommand{\elem}{\end{lem}}
\newtheorem{claim}[thm]{Claim}
\newcommand{\bclm}{\begin{claim}} \newcommand{\eclm}{\end{claim}}
\newtheorem{cor}[thm]{Corollary}
\newcommand{\bcor}{\begin{cor}} \newcommand{\ecor}{\end{cor}}
\newtheorem{conj}[thm]{Conjecture}
\newcommand{\bcnj}{\begin{conj}} \newcommand{\ecnj}{\end{conj}}
\theoremstyle{definition}
\newtheorem{defn}[thm]{Definition}
\newcommand{\bdfn}{\begin{defn}} \newcommand{\edfn}{\end{defn}}
\theoremstyle{remark}
\newtheorem{rem}[thm]{Remark}
\newcommand{\brem}{\begin{rem}} \newcommand{\erem}{\end{rem}}
\newtheorem{cnv}[thm]{Convention}
\newcommand{\bcnv}{\begin{cnv}} \newcommand{\ecnv}{\end{cnv}}
\newtheorem{exam}[thm]{Example}
\newcommand{\bexm}{\begin{exam}} \newcommand{\eexm}{\end{exam}}
\newtheorem{issue}{Issue}

\newcommand{\bpf}{\begin{proof}} \newcommand{\epf}{\end{proof}}
\newcommand{\be}{\begin{enumerate}}
\newcommand{\ee}{\end{enumerate}}
\newcommand{\bi}{\begin{itemize}}
\newcommand{\ei}{\end{itemize}}

\newcommand{\general}{\small\vfill\par\noindent\hrulefill\par
\noindent\textbf{Previous issues.} The previous issues of this
bulletin are available online at\\
\texttt{http://front.math.ucdavis.edu/search?\&t=\%22SPM+Bulletin\%22}
\\[0.1cm]
\textbf{Contributions.} Announcements, discussions, and open problems should be emailed
to \texttt{tsaban@math.biu.ac.il}\\[0.1cm]
\textbf{Subscription.}
To receive this bulletin (free) to your e-mailbox, e-mail us.
}

\newcommand{\arXiv}[5]{\subsection{#2}{#4}\par\hfill{\arx{#1}}\par\hfill\emph{#3}}

\newcommand{\AMS}[3]{\subsection{#1}\mbox{}\par\hfill{\texttt{#3}}\par\hfill\emph{#2}}

\newcommand{\arx}[1]{\texttt{http://arxiv.org/abs/#1}}
\newcommand{\url}[1]{\bq\texttt{#1}\eq}
\newcommand{\online}[1]{The paper is available online at \url{#1}}

\title[$\mathcal{SPM}$ Bulletin \textbf{\issuenumber} (\issuemonth{} \issueyear)]{%
$\mathcal{SPM}$ Bulletin\\[0.5cm]
Issue number \issuenumber: \issuemonth{} \issueyear{} \CE{}}

\begin{document}
\maketitle

\tableofcontents

\section{Editor's note}

In addition to the interesting research announcements in Section \ref{RA},
I am very pleased to announce the solution of a problem about $\gamma$-sets,
implicit since Gerlits and Nagy's 1982 paper
\emph{Some properties of $C(X)$, I},
and explicit in the 1996 paper of Just, Miller, Scheepers, and Szeptycki
\emph{The combinatorics of open covers II} and in several later papers
by these and by other authors.
Details are available in the paper announced in Section \ref{LSA} below,
and are reproduced in Section \ref{JMSSSec} below.

From a personal perspective, I am interested in problems of this sort since my
Master's thesis. In general, the question is: Assume that we take infinite
sets of natural numbers, which are rapidly thinning out in some combinatorial
sense (a scale, a tower, etc.), and then add all finite sets of naturals.
As a subspace of the Cantor space $P(\N)$, which selection hypotheses does our set satisfy?
This approach differs from the classical one, in that we do not consider the
topology during the construction. E.g., we do not take into account potential open
covers in a transfinite-inductive construction. Results of this form
were obtained by Fremllin and Miller; Just, Miller, Scheepers, and Schetycki; Scheepers;
Bartoszynski; Bartoszynski and Tsaban; and Tsaban and Zdomskyy.

The present solution, which is joint with my Master's Student Tal Orenshtein,
grew out of this series of intermediate advances, and in addition relies on the method from
Galvin and Miller's \emph{$\gamma$-sets and other singular sets of real numbers} (1984), and on 
Francis Jordan's method from \emph{There are no hereditary productive $\gamma$-spaces} (2008),
with one additional twist which makes everything fit together.
The intermediate advances which were motivated by related (but other) questions in the field of
selection principles. This is a beautiful demonstration of the importance of treating questions
in wider contexts than the ones in which they were initially posed.

Readers not interested in generalizations or in new proof methods or in weakenings of Martin's Axiom
(but which are interested in something),
may still be happy with the fact that the new result gives apparently the first proof that there are uncountable
$\gamma$-sets in all Random reals models, obtained by adding any number of Random reals
to a model of \CH{}.

\medskip

\by{Boaz Tsaban}{tsaban@math.biu.ac.il}

\hfill \texttt{http://www.cs.biu.ac.il/\~{}tsaban}

\section{$\gamma$-sets from a weak hypothesis}\label{JMSSSec}

In the paper \ref{LSA}, we construct sets of reals satisfying $\sone(\Omega,\Gamma)$,
traditionally called \emph{$\gamma$-sets}, from a weak set theoretic hypothesis.
The problem thus settled has some history, which we now survey briefly.
This involves combinatorial cardinal characteristics of the continuum \cite{BlassHBK}.
We give the necessary the definitions as we proceed.

$\gamma$-sets were introduced by Gerlits and Nagy in \cite{GN},
their most influential paper, as the third property
in a list numbered $\alpha$ through $\epsilon$. This turned out to be the most important
property in the list, and obtained its item number as it name.
One of the main results in \cite{GN} is that for Tychonoff spaces $X$,
$C(X)$ with the topology of pointwise convergence is Fr\'eche-Urysohn if, and
only if, $X$ is a $\gamma$-set.

While uncountable $\gamma$-sets exist in ZFC,\footnote{The axioms of Zermelo and Fraenkel, together with the axiom of Choice,
the ordinary axioms of mathematics.}
Borel's Conjecture (which is consistent with, but not provable within, ZFC)
implies that all metrizable $\gamma$-sets are countable.

Since we are dealing with constructions rather than general results,
we restrict attention to subsets of $\R$
(or, since the property is preserved by continuous images,
subsets of any topological space which can be embedded in $\R$).

Gerlits and Nagy proved in \cite{GN} that
Maritn's Axiom implies that all spaces of cardinality less than $\fc$ are $\gamma$-sets.
There is a simple reason for that: The \emph{critical cardinality} of a property $P$,
denoted $\non(P)$, is the minimal cardinality of a set not satisfying $P$.
Let $\binom{\Omega}{\Gamma}$ be the property: Each $\cU\in\Omega(X)$ contains a
set $\cV\in\Gamma(X)$. Gerlits and Nagy proved that $\sone(\Omega,\Gamma)=\binom{\Omega}{\Gamma}$ \cite{GN}.
Let $A\as B$ mean that $A\sm B$ is finite.
$A$ is a \emph{pseudointersection} of $\cF$ if $A\as B$ for all $B\in\cF$.
Let $\fp$ be the minimal cardinality of a family
$\cF$ of infinite subsets of $\N$ which is closed under finite intersections,
and has no pseudointersection.
Then $\non\binom{\Omega}{\Gamma}=\fp$ \cite{GM}, and Maritn's Axiom implies $\fp=\fc$ \cite{GM}.

By definition, for each property $P$, every space of cardinality smaller than $\non(P)$
satisfies $P$. Thus, the real question is whether there is $X$
of cardinality at least $\non(P)$, which satisfies $P$.
Galvin and Miller \cite{GM} proved a result of this type: $\fp=\fc$
implies that there is a $\gamma$-set of reals, of cardinality $\fp$.
Just, Miller, Scheepers and Szeptycki \cite{coc2} have improved the construction
of \cite{GM}.
We introduce their construction in a slightly more general form, that will be useful later.

Cantor's space $\Cantor$ is equipped with the Tychonoff product topology,
and $P(\N)$ is identified with $\Cantor$ using characteristic functions.
This defines the topology of $P(\N)$. The partition $P(\N)=\roth\cup\Fin$, into
the infinite and the finite sets, respectively, is useful here.

For $f,g\in\NN$, let $f\le^* g$ if $f(n)\le g(n)$ for all but finitely many
$n$. $\fb$ is the minimal cardinality of a $\le^*$-unbounded subset of $\NN$.
A set $B\sbst\roth$ is \emph{unbounded} if the set of all increasing enumerations
of elements of $B$ is unbounded in $\NN$, with respect to $\le^*$.

\bdfn
A\emph{tower} of cardinality $\kappa$ is a set $T\sbst\roth$ which can be enumerated
bijectively as $\{x_\alpha : \alpha<\kappa\}$,
such that for all $\alpha<\beta<\kappa$, $x_\beta\as x_\alpha$.

An \emph{unbounded tower} of cardinality $\kappa$ is an unbounded set $T\sbst\roth$ which is
a tower of cardinality $\kappa$. (Necessarily, $\kappa\ge\fb$.)
\edfn

Let $\ft$ be the minimal cardinality of a tower which has no pseudointersection.
Rothberger proved that $\ft\le\fb$ \cite{BlassHBK}.
$\ft=\fb$ if, and only if, there is an unbounded tower of cardinality $\ft$.

Just, Miller, Scheepers and Szeptycki \cite{coc2} proved that if
$T$ is an unbounded tower of cardinality $\aleph_1$,
then $T\cup\Fin$ satisfies $\sone(\Omega,\Omega)$, as well as
a property, which was later proved by Scheepers \cite{wqn} to be equivalent to $\sone(\Gamma,\Gamma)$ .
In Problem 7 of \cite{coc2}, we are asked the following.

\bprb[Just-Miller-Scheepers-Szeptycki \cite{coc2}]\label{JMSSProb}
Assume that $T\sbst\roth$ is an unbounded tower of cardinality $\aleph_1$ (so that $\aleph_1=\fb$).
Is $T\cup\Fin$ a $\gamma$-set, i.e., satisfies $\sone(\Omega,\Gamma)$?
\eprb

Scheepers proves in \cite{alpha_i} that for each unbounded tower $T$ of cardinality $\ft=\fb$,
$T\sbst\roth$ satisfies $\sone(\Gamma,\Gamma)$.

Miller \cite{MillerBC} proves that in
the Hechler model, there are no uncountable $\gamma$-sets. In this model, $\aleph_1=\fp=\ft<\fb$,
and thus $\aleph_1=\ft$ does not suffice to have an uncountable $\gamma$-set of reals.
At the end of \cite{MillerBC} and in its appendix, Miller proves that $\lozenge(\fb)$,
a property strictly stronger than $\aleph_1=\fb$, implies that there is an uncountable $\gamma$-set of reals.\footnote{$\lozenge(\fb)$ is defined in Dzamonja-Hrusak-Moore \cite{DHM}.}
He concludes that it is still open whether $\fb = \aleph_1$ is enough to construct an uncountable $\gamma$-set.

We show that the answer is positive, and indeed also answer a question of Gruenhage and Szeptycki \cite{FUfin}:
A classical problem of Malykhin asks whether there is a countable Fr\'echet-Urysohn topological group which is not
metrizable. Gruenhage and Szeptycki prove that $F\sbst\NN$ is a $\gamma$-set if, and only if, a certain construction
associated to $F$ provides a positive answer to Malykhin's Problem \cite{FUfin}.
They define a generalization of $\gamma$-set, called \emph{weak $\gamma$-set}, and combine their results with results of Nyikos
to prove that $\fp=\fb$ implies that there is a weak $\gamma$-set in $\NN$ \cite[Corollary 10]{FUfin}.
They write: ``The relationship between $\gamma$-sets and weak $\gamma$-sets is not known.
Perhaps $\fb = \fp$ implies the existence of a $\gamma$-set.'' Our solution confirms their conjecture.

$\fp\le\ft\le\fb$, and in all known models of set theory, $\fp=\ft$.
Our theorem reproduces Galvin and Miller's Theorem when $\fp=\fc$ \cite{GM}, but
gives additional information: Even if the
possible open covers are not considered during the construction in \cite{GM},
the resulting set is still a $\gamma$-set.

\bthm\label{JMSSSol}
For each unbounded tower $T$ of cardinality $\fp$ in $\roth$,
$T\cup\Fin$ satisfies $\sone(\Omega,\Gamma)$.
\ethm

Zdomskyy points out that our proof actually shows that a wider family of sets are $\gamma$-sets.
For example, if we start with $T$ an unbounded tower of cardinality $\fp$, and thin out its elements
arbitrarily, $T\cup\Fin$ remains a $\gamma$-set. This may be useful for constructions of examples
with additional properties, since this way, each element of $T$ may be chosen arbitrarily from a certain perfect set.

In particular, we have that in each model of ZFC where $\fp=\fb$, there are $\gamma$-sets of cardinality $\fp$.

\bcor\label{ec}
In each of the Cohen, Random, Sacks, and Miller models of ZFC, there are $\gamma$-sets of reals
with cardinality $\fp$.\qed
\ecor

As discussed above, there are no uncountable $\gamma$-sets in the Hechler model \cite{MillerBC}.
Since the Laver and Mathias models satisfy Borel's Conjecture, there are no
uncountable $\gamma$-sets in these models, too.

Earlier, Corollary \ref{ec} was shown for the Sacks model by Ciesielski, Mill\'an, and Pawlikowski
in \cite{CPAgamma},
and for the Cohen and Miller models by Miller \cite{MillerBC}, using specialized
arguments. It seems that the result, that there are uncountable $\gamma$-sets
in the Random reals model (constructed by extending a model of \CH{}), is new.

\nby{Tal Orenshtein and Boaz Tsaban}

\section{Research announcements}\label{RA}

\AMS{Ultrafilters with property $(s)$}{Arnold W. Miller}{http://www.ams.org/journal-getitem?pii=S0002-9939-09-09919-5}

\AMS{On a converse to Banach's Fixed Point Theorem}
{Marton Elekes}
{http://www.ams.org/journal-getitem?pii=S0002-9939-09-09904-3}

\AMS{Analytic groups and pushing small sets apart}
{Jan van Mill}
{http://www.ams.org/journal-getitem?pii=S0002-9947-09-04665-0}

\arXiv{0905.3754}
{Club-guessing, stationary reflection, and coloring theorems}
{Todd Eisworth}
{We obtain strong coloring theorems at successors of singular cardinals from
failures of certain instances of simultaneous reflection of stationary sets.
Along the way, we establish new results in club-guessing and in the general
theory of ideals.}

\arXiv{0905.3913}
{More on the pressing down game}
{Jakob Kellner, Saharon Shelah}
{We investigate the pressing down game and its relation to the Banach Mazur
game. In particular we show: Consistently relative to a supercompact, there is
a nowhere precipitous normal ideal $I$ on $\aleph_2$ such that player nonempty
wins the pressing down game of length $\aleph_1$ on $I$ even if player empty
starts. For the proof, we construct a forcing notion to force the following:
There is normal, nowhere precipitous ideal $I$ on a supercompact $\kappa$ such
that for every $I$-positive $A$ there is a normal ultrafilter containing $A$
and extending the dual of $I$.}

\arXiv{0905.3588}
{A note on discrete sets}
{Santi Spadaro}
{We give several partial positive answers to a question of Juhasz and
Szentmiklossy regarding the minimum number of discrete sets required to cover a
compact space. We study the relationship between the size of discrete sets,
free sequences and their closures with the cardinality of a Hausdorff space,
improving known results in the literature.}

\AMS{Antidiamond principles and topological applications}
{Todd Eisworth and Peter Nyikos}
{http://www.ams.org/journal-getitem?pii=S0002-9947-09-04705-9}

\arXiv{0907.3771}
{Partitions and indivisibility properties of countable dimensional vector spaces}
{C. Laflamme, L. Nguyen Van The, M. Pouzet, N. Sauer}
{We investigate infinite versions of vector and affine space partition
results, and thus obtain examples and a counterexample for a partition problem
for relational structures. In particular we provide two (related) examples of
an age indivisible relational structure which is not weakly indivisible.}

\arXiv{0907.4941}
{Group-valued continuous functions with the topology of pointwise convergence}
{Dmitri Shakhmatov, Jan Spevak}
{Let $G$ be a topological group with the identity element $e$. Given
a space $X$, we denote by $C_p{X}{G}$ the group of all continuous
functions from $X$ to $G$ endowed with the topology of pointwise
convergence, and we say that $X$ is: (a) {\em $G$-regular\/} if, for
each closed set $F\subseteq X$ and every point $x\in X\setminus F$,
there exist $f\in C_p{X}{G}$ and $g\in G\setminus\{e\}$ such that
$f(x)=g$ and $f(F)\subseteq\{e\}$; (b) {\em $G^\star$-regular\/}
provided that there exists $g\in G\setminus\{e\}$ such that, for
each closed set $F\subseteq X$ and every point $x\in X\setminus F$,
one can find $f\in C_p{X}{G}$ with $f(x)=g$ and $f(F)\subseteq\{e\}$.
 Spaces $X$ and $Y$ are {\em $G$-equivalent\/} provided that the
topological groups $C_p{X}{G}$ and $C_p{Y}{G}$ are topologically
isomorphic.

We investigate which topological properties are preserved by
$G$-equivalence, with a special emphasis being placed on
characterizing topological properties of $X$ in terms of those of
$C_p{X}{G}$. Since $\mathbb{R}$-equivalence coincides with
$l$-equivalence, this line of research ``includes'' major topics of
the classical $C_p$-theory
 of Arhangel'ski\u{\i} as a particular case (when $G=\mathbb{R}$).

We introduce a new class of TAP groups that contains all groups
having no small subgroups (NSS groups). We prove that: (i) for a
given NSS group $G$, a $G$-regular space $X$ is pseudocompact if and
only if $C_p{X}{G}$ is TAP, and (ii) for a metrizable NSS group
$G$, a
 $G^\star$-regular space $X$ is compact if and only if $C_p{X}{G}$ is
a TAP group of countable tightness. In particular, a Tychonoff
space $X$ is pseudocompact (compact) if and only if
$C_p{X}{\mathbb{R}}$ is a TAP group (of countable tightness).
Demonstrating the limits of the result in (i), we give an example of
a precompact TAP group $G$ and a $G$-regular countably compact space
$X$ such that $C_p{X}{G}$ is not TAP.

We show that Tychonoff spaces $X$ and $Y$ are
$\mathbb{T}$-equivalent if and only if their free precompact Abelian
groups are topologically isomorphic, where $\mathbb{T}$ stays for
the quotient group $\mathbb{R}/\mathbb{Z}$. As a corollary, we
obtain that $\mathbb{T}$-equivalence implies $G$-equivalence for
every Abelian precompact group $G$. We
establish that $\mathbb{T}$-equivalence preserves the following
topological properties: compactness, pseudocompactness,
$\sigma$-compactness, the property of being a Lindel\"of
$\Sigma$-space, the property of being a compact metrizable space,
the (finite) number of connected components, connectedness, total
disconnectedness. An example of $\mathbb{R}$-equivalent (that is,
$l$-equivalent) spaces that are not $\mathbb{T}$-equivalent is
constructed.}

\arXiv{0907.4126}
{Stationary and convergent strategies in Choquet games}
{Fran\c{c}ois G. Dorais and Carl Mummert}
{If POINT has a winning strategy against EMPTY in the
  Choquet game on a space, the space is said to be a
  Choquet space.  Such a winning strategy allows
  POINT to consider the entire finite history of previous
  moves before making each new move; a stationary strategy only
  permits POINT to consider the previous move by EMPTY.
  We show that POINT has a stationary winning strategy for
  every second countable $T_1$ Choquet space. More generally,
  POINT has a stationary winning strategy for any $T_1$ Choquet
  space with an open-finite basis.

  We also study convergent strategies for the Choquet game,
  proving the following results.

A $T_1$ space $X$ is the open
  image of a complete metric space if and only if POINT
  has a convergent winning strategy in the Choquet game
  on~$X$.

A $T_1$ space $X$ is the
  compact open image of a metric space if and only if $X$ is
  metacompact and POINT has a stationary convergent strategy in
  the Choquet game on~$X$.

A $T_1$ space $X$ is the
  compact open image of a complete metric space if and only if $X$ is
  metacompact and POINT has a stationary convergent winning strategy in
  the Choquet game on~$X$.
}

\arXiv{0906.5136\label{LSA}}
{Linear $\sigma$-additivity and some applications}
{Tal Orenshtein and Boaz Tsaban}
{We show that countable
increasing unions preserve a large family of well-studied covering
properties, which are not necessarily $\sigma$-additive.
Using this, together with infinite-combinatorial methods and simple forcing theoretic
methods, we explain several phenomena, settle problems of
Just, Miller, Scheepers and Szeptycki;
Gruenhage and Szeptycki;
Tsaban and Zdomskyy;
and Tsaban,
and construct topological groups with very strong combinatorial properties.
(See also Sections 1 and 2 above.)
}

\ed
\begin{thebibliography}{99}

\bibitem{BlassHBK}
A. Blass,
\emph{Combinatorial cardinal characteristics of the continuum},
in: \textbf{Handbook of Set Theory} (M.\ Foreman, A.\ Kanamori, and M.\ Magidor, eds.),
Kluwer Academic Publishers, Dordrecht, to appear.
\texttt{http://www.math.lsa.umich.edu/\~{}ablass/hbk.pdf}

\bibitem{CPAgamma}
K. Ciesielski, A. Mill\'an, and J. Pawlikowski,
\emph{Uncountable $\gamma$-sets under axiom $\mathrm{CPA}^\mathrm{game} \mathrm{cube}$},
Fundamenta Mathematicae \textbf{176} (2003),  143--155.

\bibitem{DHM}
M. D\v{z}amonja, M. Hru\'sak, and J. Moore,
\emph{Parametrized $\lozenge$ principles},
Transactions of the American Mathematical Society \textbf{356} (2004), 2281--2306.

\bibitem{GM}
F.\ Galvin and A. Miller,
\emph{$\gamma$-sets and other singular sets of real numbers},
Topology and its Applications \textbf{17} (1984), 145--155.

\bibitem{GN}
J.\ Gerlits and Zs.\ Nagy,
\emph{Some properties of $C(X)$, I},
Topology and its Applications \textbf{14} (1982), 151--161.

\bibitem{FUfin}
Gary Gruenhage and Paul Szeptycki,
\emph{Fr\'echet–Urysohn for finite sets},
Topology and its Applications \textbf{151} (2005) 238--259.

\bibitem{coc2}
W. Just, A. Miller, M. Scheepers, and P. Szeptycki,
\emph{The combinatorics of open covers II},
Topology and its Applications \textbf{73} (1996), 241--266.

\bibitem{MillerBC}
A. Miller,
\emph{The $\gamma$-Borel conjecture},
Archive for Mathematical Logic \textbf{44} (2005), 425--434.
\arx{math.LO/0312308}

\bibitem{alpha_i}
M.\ Scheepers,
\emph{$C_p(X)$ and Arhangel'ski\u{\i}'s $\alpha_i$ spaces},
Topology and its Applications \textbf{89} (1998), 265--275.

\bibitem{wqn}
M.\ Scheepers,
\emph{Sequential convergence in ${\sf C}_p(X)$ and a covering property},
East-West Journal of Mathematics \textbf{1} (1999),
207--214.

\end{thebibliography}
